\documentclass[a4paper]{amsart}

\pdfoutput=1

\usepackage{amsthm, amssymb, amsmath, amsfonts, mathrsfs}

\usepackage{microtype}

\usepackage[pagebackref,colorlinks=true,pdfpagemode=none,urlcolor=blue,
linkcolor=blue,citecolor=blue]{hyperref}

\usepackage{relsize}

\usepackage{amsmath,amsfonts,amssymb,amsthm}
\usepackage{amsthm, amssymb, amsmath, amsfonts, mathrsfs}

\usepackage{mathrsfs}
\usepackage{MnSymbol}
\usepackage{scalerel} 

\usepackage{color}
\usepackage{accents}

\usepackage[normalem]{ulem}
\usepackage{bbm}

\usepackage{cleveref}

\usepackage{mathtools}

\definecolor{labelkey}{gray}{.8}
\definecolor{refkey}{gray}{.8}

\definecolor{darkred}{rgb}{0.9,0.1,0.1}
\definecolor{darkgreen}{rgb}{0,0.5,0}

\newtheorem{theorem}{Theorem}[section]

\theoremstyle{remark}
\newtheorem{remark}[theorem]{Remark}

\numberwithin{equation}{section}
\newcommand{\R}{\mathbb{R}}

\newcommand{\bbR}{\mathbb{R}}

\newcommand{\Z}{\mathbb{Z}}

\newcommand{\PP}{\mathbf{P}}
\newcommand{\E}{\mathbb{E}}

\newcommand{\F}{\mathcal{F}}

\newcommand{\cN}{\mathcal{N}}

\newcommand{\EE}{\mathbf{E}}

\newcommand{\cP}{\mathcal{P}}

\begin{document}
\title{Noise sensitivity for stochastic heat and Schr\"odinger equation}

\author{Yu Gu, Tomasz Komorowski}

\address[Yu Gu]{Department of Mathematics, University of Maryland, College Park, MD 20742, USA. }
\email{yugull05@gmail.com}

\address[Tomasz Komorowski]{Institute of Mathematics, Polish Academy of Sciences, ul.
Śniadeckich 8, 00-636 Warsaw, Poland. }
\email{tkomorowski@impan.pl}

\maketitle

\begin{abstract}

In this note, we consider the  stochastic heat and Schr\"odinger equation, and show that, at time $t$, the onset of the chaos occurs on the scale of $1/t$, and the Fourier spectrum  of the solution is asymptotically Gaussian after centering and rescaling. 
\medskip


\end{abstract}

\section{Introduction and main result}

\subsection{Fourier spectrum of a random variable with respect to a
  Wiener chaos} Suppose that $\big(V(t,x)\big)_{(t,x)\in\bbR^{1+d}}$ is a centered
stationary Gaussian random field built on some probability space
$(\Omega,\F,\PP)$. Its covariance is of the form 
\[
\EE \big[V(t,x)V(s,y)\big]=\delta(t-s)R(x-y),
\]
 where $R(\cdot)$ is the spatial covariance function that is either
 smooth and fast decaying in $x$ (colored noise), or $R(x)=\delta(x)$
 (white noise).
 Suppose that a random variable we
 consider  $X_t$ where $t$ is the time variable which will become
 large, and assume $X_t$ is measurable with respect to $\F_t$ which is
 the $\sigma-$algebra generated by the noise up to time $t$. If
 $X_t\in L^2(\Omega)$, it can be expanded into the Wiener chaos
 associated with $V$, see an introduction in \cite[Chapter
 1]{nualart}. In this particular case, when   $V$ is white in the time, it can be written as $X_t=\sum_{n\geq0} I_n (f_{t,n})$ with 
 \begin{equation}
   \label{010402-25}
\begin{aligned}
 I_n(f_{t,n}):=\int_{[0,t]_<^n} \int_{\R^{nd}}f_{t,n}(s_{1,n};y_{1,n})V(s_1,y_1)\ldots V(s_n,y_n) dy_{1,n}ds_{1,n},
\end{aligned}
\end{equation}
where $[0,t]_<^n=\{(s_1,\ldots,s_n)\in \R_+^n: 0<s_n<\ldots<s_1<t\}$,
\begin{align*}
 & s_{1,n}=(s_1,\ldots,s_n), \quad ds_{1,n}=ds_1\ldots ds_n,\\
&y_{1,n}=(y_1,\ldots,y_n),\quad dy_{1,n}=dy_1\ldots dy_n.
\end{align*}
The  stochastic integral appearing in \eqref{010402-25} is interpreted
as the iterated It\^o integral. With the above representation, the
second moment of $X_t$ can be written as 
 \begin{equation}\label{e.defctn}
   \begin{aligned}
     \EE X_t^2&=\sum_{n\geq0} c_{t,n}^2,\quad \mbox{ with}\\
 c_{t,n}^2&= \int_{[0,t]_<^n}\int_{\R^{2nd}}f_{t,n}(s_{1,n}
 ;y_{1,n})f_{t,n}^*(s_{1,n} ;\tilde y_{1,n})&
 \\
 &\qquad \qquad \qquad \times\prod_{j=1}^n  R(y_j-\tilde{y}_j) dy_{1,n}d\tilde{y}_{1,n}ds_{1,n}.
 \end{aligned}
 \end{equation}
 In other words, $c_{t,n}^2$ is the contribution of the $n-$th Wiener
 chaos to $\EE X_t^2$. Following noise sensitivity literature (see
 the book \cite{GS14}), we define the Fourier spectrum of the random
 variable $X_t$ as the  probability distribution on non-negative
 integers $\Z_{\geq0}$: 
 \begin{equation}
   \label{030602-25}
 \PP(\mathcal{N}_t=n)=\frac{c_{t,n}^2}{\EE X_t^2}, \quad\quad n\geq0.
 \end{equation}
In our presentation, we shall also call any random variable
$\mathcal{N}_t$, having the aforementioned distribution, the {\em Fourier
spectrum} of $X_t$. We are particularly interested in the asymptotics
of $\mathcal{N}_t$, as
$t\to\infty$, which describes the   composition of
the second moment, in terms of  the Wiener chaos of different orders.


\subsection{Random Schr\"odinger and stochastic heat equations}

 As  concrete examples  we consider two stochastic linear equations
 with  multiplicative noises:
\begin{equation}\label{e.schrshe}
\begin{array}{ll}
i\partial_t\phi+\tfrac12\Delta\phi+V(t,x)\phi=0, &\phi(0,x)=\phi_0(x),\\
                                                 & \\
  \partial_t Z=\tfrac12\Delta Z+V(t,x)Z, &
Z(0,x)=Z_0(x),\,(t,x)\in(0,+\infty)\times \R^d.
\end{array}
\end{equation}
They are referred to as the random Schr\"odinger and stochastic heat
equation, respectively.
We will specify later on the assumptions made about the spatial
dimension $d$, the spatial covariance function of the noise $R(\cdot)$
and the initial data $\phi_0(\cdot)$, $Z_0(\cdot)$.
For the moment we suppose that they suffice to establish the existence
of a unique random field solution to a particular equation of
\eqref{e.schrshe}, and the solution evaluated at any spacetime point
is a random variable that belongs to  $L^2(\Omega)$. We are interested
in the chaotic behaviors of such $L^2(\Omega)$ valued random elements.
  
Consider two cases.
\begin{itemize}
\item[case 1)] We have  $d\geq 1$, both the covariance function $R(\cdot)$ and the initial
  data $\phi_0(\cdot)$ belong to  the Schwartz class
  $\mathcal{S}(\R^d)$,
  and 
  \begin{equation}
    \label{e.defXt1}
    X_t= \hat{\phi}(t,0).
  \end{equation}
  Here and throughout the paper the Fourier transform is defined as
  $$
  \hat{f}(\xi)=\int_{\R^d} f(x)e^{-i\xi\cdot x}dx,\quad \mbox{for any }f\in L^1(\R^d).
  $$
 
\item[case 2)] the spatial dimension $d=1$, $R(\cdot)=\delta(\cdot)$,
  $Z_0(\cdot)\equiv1$ and 
  \begin{equation}
    \label{e.defXt}
    X_t=  Z(t,0).
\end{equation}
  \end{itemize}
  The assumptions on $R(\cdot), \phi_0(\cdot),Z_0(\cdot)$ can be
  greatly relaxed. Our choice is motivated by clarity of   the presentation. 
  
  Here is the main result:
  \begin{theorem}\label{t.mainth}
  Let $\mathcal{N}_t$ be the Fourier spectrum associated with $X_t$
  defined in either \eqref{e.defXt1} or \eqref{e.defXt}. In both cases, there exist $\sigma,\mu>0$ such that
  \begin{equation}\label{e.clt}
  \frac{\mathcal{N}_t-\mu t}{\sigma\sqrt{t}}\Rightarrow N(0,1),\quad\mbox{ as }t\to\infty
  \end{equation}
  in law. Here $N(0,1)$ is the standard normal law (i.e. mean
  zero and variance $1$).
  \end{theorem}

  \subsection{Noise sensitivity} Theorem~\ref{t.mainth} shows that $\mathcal{N}_t$ concentrates on the scale of order $t$, meaning the dominant contribution to the second moment of $X_t$ comes from Wiener chaos of order  $t$. This proves the noise sensitivity. More precisely, consider the following perturbation of $V$: for any $s>0$, define the perturbation of strength $s$ by 
\begin{equation}\label{e.defVs}
V_s(t,x)=e^{-s} V(t,x)+\sqrt{1-e^{-2s}}\tilde{V}(t,x),
\end{equation}
where $\tilde{V}(\cdot,\cdot)$ is an independent copy of
$V(\cdot,\cdot)$. Let $X_t(s)$ be   the analogue of $X_t$ obtained
from \eqref{e.schrshe} with the noise $V$ replaced by $V_s$. In other
words, we consider the stochastic Schr\"odinger and heat equation
driven by $V_s$ and $X_t(s)$ is the corresponding $\hat{\phi}(t,0)$
and $Z(t,0)$. The quantity of interest is the correlation between
$X_t(s)$ and $X_t(0)$, and in particular the size of the  perturbation
required so that we observe a strict decorrelation. The correlation
can be expressed in terms of the coefficients $c_{t,n}^2$ as follows:
\[
\begin{aligned}
\mathrm{Cor}[X_t(0),X_t(s)]=&\frac{\EE [X_t(0)X_t^*(s)]-\EE [X_t(0)\EE X_t^*(s)]}{\EE|X_t(0)-\EE X_t(0)|^2}\\
&=\frac{\sum_{n\geq1} e^{-ns} c_{t,n}^2}{\sum_{n\geq1} c_{t,n}^2}=\EE[e^{-s\mathcal{N}_t}\,|\, \mathcal{N}_t>0]=\frac{\EE e^{-s\mathcal{N}_t}-\PP(\mathcal{N}_t=0)}{1-\PP(\mathcal{N}_t=0)}.
\end{aligned}
\]
If we consider perturbation corresponding to $s\sim t^{-\alpha}$ for some $\alpha>0$, a direct corollary of Theorem~\ref{t.mainth} is that
\[
\mathrm{Cor}[X_t(0),X_t(s)]\to \left\{
\begin{array}{ll} 
1 & \mbox{ if } \alpha>1,\\
0 & \mbox{ if } \alpha<1,
\end{array}
\right.
\]
as $t\to\infty$.
In other words, the onset of chaos is on the scale of $s\sim
t^{-1}$. Here the symbol $s \sim f(t)$, for
a positive valued function $f(\cdot)$ means that there exists a
constant $c\in(0,1)$ such that $c f(t)\le s\le c^{-1}f(t)$ for
sufficiently large $t$.

The equations we study in this note are simple enough that the
expansion coefficients are explicit,  allowing for direct
calculations. Still, we believe it is interesting to present the
phenomenon   quantitatively here. In both cases the random variable
$\mathcal{N}_t$ turns out to be  close to a Poisson distribution with
intensity of order $t$, so \eqref{e.clt} may be interpreted as a
central limit theorem for an approximate Poisson distributions with a
large intensity parameter. In fact, the mechanism underlying the chaotic behavior for \eqref{e.schrshe} is analogous to that of a geometric Brownian motion: for $X_t=\exp(B_t-\tfrac12t)$ where $B$ is a standard Brownian motion, $\mathcal{N}_t$ follows a Poisson distribution with intensity $t$.

It would be   interesting to study chaotic behaviors of solutions to nonlinear SPDEs and consider random variables such as $X_t=\log Z(t,0)$, which can be interpreted as the free energy of a continuum directed random polymer. A plausible conjecture is $N_t\sim t^{1/3}$, but proving this remains a challenging open problem. For relevant references, one may consult the book \cite{chatterjee}, and the recent work \cite{GH24}, which studied a zero-temperature model of Brownian last passage percolation and proved a similar result for the overlap of geodesics.

\subsection*{Acknowledgement}   Y. G. was partially supported by the NSF
through DMS-2203014.  T. K. acknowledges the support of the NCN grant 2020/37/B/ST1/00426.

\section{It\^o-Schr\"odinger equation}
\label{s.schr}

In this section we prove Theorem~\ref{t.mainth} for the
It\^o-Schr\"odinger equation
\begin{equation}
  \label{ito-s}
  i\partial_t\phi+\tfrac12\Delta\phi+V(t,x)\circ \phi=0.
\end{equation}
Here the product $\circ$ between the noise $V(\cdot,\cdot)$ and $\phi$ is
 interpreted in the Stratonovitch sense so that the equation
preserves the $L^2(\R^d)$ norm. 

%
The equation \eqref{ito-s} can be written in the Fourier domain and in It\^o form as 
 \[
 i\partial_t\hat{\phi}(t,\xi)-\tfrac12|\xi|^2\hat{\phi}(t,\xi)+\tfrac{i}{2}R(0)\hat{\phi}(t,\xi)-\int_{\R^d}\hat{\phi}(t,\xi-p) \frac{\hat{V}(t,p) dp}{(2\pi)^d} =0,
 \]
 see \cite[Section 5.1]{GK21}. Here $\hat{V}(t,p)$ is the (complex-valued) Gaussian noise with the correlation
\begin{equation}
\label{corel}
\begin{aligned}
&\EE\left[\hat{V}(t,p)
 \hat{V}^*(s,q)\right]=(2\pi)^d\hat R(p)\delta(t-s)\delta(p-q),\\
&
\hat{V}^*(t,p)=\hat{V}(t,-p).
\end{aligned}
\end{equation} 
Define the {\em compensated wave function}
\begin{equation}
  \label{010602-25}
\psi(t,\xi)=\hat{\phi}(t,\xi)\exp\left\{\frac12\big(R(0)+i|\xi|^2\big)t\right\}.
\end{equation}
It satisfies 
\[
\begin{aligned}
&\partial_t \psi(t,\xi)=\int_{\R^d}\psi(t,\xi-p)\exp\left\{\frac{i}{2}( |\xi|^2-|\xi-p|^2)t\right\}\frac{\hat{V}(t,p) dp}{i(2\pi)^d},\\
&\psi(0,\xi)=\hat{\phi}_0(\xi).
\end{aligned}
\]
Through an iteration,  $\psi(t,\xi)$ can be written as the Wiener chaos expansion $\psi(t,\xi)=\sum_{n\geq0} \psi_n(t,\xi)$ with $\psi_0(t,\xi)=\hat{\phi}_0(\xi)$ and 
\begin{equation}
  \label{020602-25}
  \begin{split}
\psi_n(t,\xi)=\int_{[0,t]_<^n}\int_{\R^{nd}}&\hat{\phi}_0(\xi-p_1-\ldots-p_n)
\\
&
\times\prod_{j=1}^n
 \Bigg[\exp\left\{\frac{i}{2}(|\xi-p_0\ldots-p_{j-1}|^2-|\xi-p_0\ldots-p_j|^2)s_j\right\}\frac{\hat{V}(s_j,p_j)}{i(2\pi)^d}\Bigg]dp_{1,n}ds_{1,n},
\end{split}
\end{equation}
with the convention $p_0=0$. Therefore (cf \eqref{010602-25}) $\hat{\phi}$ admits the  chaos expansion 
 \[
 \hat{\phi}(t,\xi)=\psi(t,\xi)e^{-\tfrac12(R(0)+i|\xi|^2)t}=\sum_{n\geq0} \psi_n(t,\xi)e^{-\tfrac12(R(0)+i|\xi|^2)t}.
 \] 
From the explicit expression \eqref{020602-25}, we  obtain
\begin{equation}\label{e.251}
\begin{aligned}
\EE |\psi_n(t,\xi)|^2=\tfrac{t^n}{n!}\int_{\R^{nd}}\hat{\Phi}_0 (\xi-p_1-\ldots-p_n)\prod_{j=1}^n \frac{\hat{R}(p_j)}{(2\pi)^d} dp_{1,n},
\end{aligned}
\end{equation}
where
$
 \hat{\Phi}_0(\xi)=|\hat{\phi}_0(\xi)|^2. 
 $
 
 Recall, see \eqref{030602-25}, that  the Fourier spectrum associated with
 the  random variable $\hat{\phi}(t,0)$ is  given
 by  \[
 \PP(\mathcal{N}_t=n)=\frac{c_{t,n}^2}{\sum_{j\geq0} c_{t,j}^2},
 \]
 with 
 \[
 \begin{aligned}
 &c_{t,0}^2=\hat{\Phi}_0(0) e^{-R(0)t},\\
 &c_{t,n}^2=\frac{t^n}{n!}e^{-R(0)t}\int_{\R^{nd}}  \hat{\Phi}_0 (-p_1-\ldots-p_n)  \prod_{j=1}^n \frac{\hat{R}(p_j) }{(2\pi)^d}  dp_{1,n}, \quad\quad n\geq1.
 \end{aligned}
\]

 To prove the convergence in distribution of the rescaled $\mathcal{N}_t$, we write:
 \[
 \int_{\R^{nd}} \prod_{j=1}^n \hat{\Phi}_0 (-p_1-\ldots-p_n) \frac{\hat{R}(p_j)}{(2\pi)^d} dp_{1,n}= \int_{\R^d} R(y)^n\Phi_0(y)dy,
 \]
 where $\Phi_0$ is the inverse Fourier transform of $\hat\Phi_0$. Thus,
 \[
 c_{t,n}^2=\frac{t^n}{n!}e^{-R(0)t}\int_{\R^d} R(y)^n \Phi(y)_0dy, \quad\quad n\geq0.
 \]
 For any $\theta\in\R$, one can compute therefore
 \[
 \begin{aligned}
 \EE e^{i\theta \mathcal{N}_t}=\frac{\sum_{n=0}^\infty e^{i\theta n}c_{t,n}^2}{\sum_{n=0}^\infty c_{t,n}^2}=\frac{\int_{\R^d} \exp\left\{t(R(y)e^{i\theta}-R(0))\right\}\Phi_0(y)dy}{\int_{\R^d} \exp\left\{t(R(y)-R(0)) \right\}\Phi_0(y)dy }.
 \end{aligned}
 \]
 In consequence
 \[
 \begin{aligned}
 &\EE \exp\left\{i\theta \frac{\mathcal{N}_t-R(0)t}{\sqrt{t}}\right\}\\
 &=\frac{\int_{\R^d} \exp\left\{t(R(y)e^{i\theta/\sqrt{t}}-R(0))\right\}\Phi_0(y)dy}{\int_{\R^d} \exp\left\{t(R(y)-R(0)) \right\}\Phi_0(y)dy }\exp\left\{-i\theta R(0)\sqrt{t}\right\}\\
 &=\frac{\int_{\R^d} \exp\left\{t(R(y/\sqrt{t})e^{i\theta/\sqrt{t}}-R(0))\right\}\Phi_0(y/\sqrt{t})dy}{\int_{\R^d} \exp\left\{t(R(y/\sqrt{t})-R(0)) \right\}\Phi_0(y/\sqrt{t})dy }\exp\left\{-i\theta R(0)\sqrt{t}\right\}.
 \end{aligned}
 \]
 Since
 $\nabla R(0)=0$ (due to the parity of $R(\cdot)$), by the Taylor expansion,   the denominator converges to
 \[
 \int_{\R^d} \exp\left\{t(R(y/\sqrt{t})-R(0))
 \right\}\Phi_0(y/\sqrt{t})dy \to \Phi_0(0)\int_{\R^d}\exp\left\{\tfrac12y\cdot
 \nabla^2R(0) y\right\} dy,\quad t\to+\infty,
 \]
 where $\nabla^2R(0)$ is the Hessian
 matrix, which we assume is
 strictly negative definite.  We have $\Phi_0(0)=\tfrac{1}{(2\pi)^d}\int_{\R^d}|\hat{\phi}_0(\xi)|^2d\xi>0$.
 Similarly,   the numerator 
 \[
 \begin{aligned}
 &\int_{\R^d}
 \exp\left\{t(R(y/\sqrt{t})e^{i\theta/\sqrt{t}}-R(0))\right\}\Phi_0(y/\sqrt{t})dy\exp\left\{-i\theta
   R(0)\sqrt{t}\right\}\\
 &
 =\int_{\R^d}
   \exp\left\{tR(y/\sqrt{t})\big(e^{i\theta/\sqrt{t}}-1\big)\right\}\\
   &
   \times \exp\left\{t(R(y/\sqrt{t})-R(0))
 \right\}\Phi_0(y/\sqrt{t})dy\exp\left\{-i\theta R(0)\sqrt{t}\right\}\\
 &\to \exp\left\{-\tfrac12R(0)\theta^2\right\}\Phi_0(0) \int_{\R^d}\exp\left\{\tfrac12y\cdot R''(0) y\right\} dy.
 \end{aligned}
 \]
 We have shown therefore 
 \[
 \EE \exp\left\{i\theta \tfrac{\mathcal{N}_t-R(0)t}{\sqrt{t}})\to \exp(-\tfrac12R(0)\theta^2\right\},
 \]
 which completes the proof of Theorem \ref{t.mainth} in the case of
 the solution of a random Schr\"odinger equation.

 \begin{remark}
 It is known that the second moment $w(t,\xi):=\EE |\hat{\phi}(t,\xi)|^2$ solves the  kinetic equation 
\begin{equation}\label{e.kinetic}
\partial_t w(t,\xi)=\int_{\R^d} [w(t,\xi-p)-w(t,\xi)]\frac{\hat{R}(p) }{(2\pi)^d}dp, \quad\quad w(0,\xi)=|\hat{\phi}_0(\xi)|^2.
\end{equation}
 The underlying probabilistic model is a compound Poisson process with
intensity
$R(0)$ and the jump size distribution
is given by the density function $\bar R(p)=\frac{\hat{R}(p) }{(2\pi)^dR(0)}$.
Under the diffusive scaling, the solution of \eqref{e.kinetic}
converges to the solution of the heat equation with the diffusivity
matrix $(2\pi)^{-d}\int_{\R^d}p_jp_k \hat{R}(p)dp$, $j,k=1,\ldots,d$. The respective scaling limit
of the underlying Poisson 
process is a Brownian
motion. One should view then $\mathcal{N}_t$   as the number of  
scatterings (or jumps) of the Poisson process, which explains the aforementioned central limit theorem. 
  \end{remark}

\section{Stochastic heat equation}

In this section we prove Theorem~\ref{t.mainth} for the case of $X_t=Z(t,0)$, with $Z(t,x)$ solving the stochastic heat equation. For the proof, we consider the more general case when there is an inverse temperature parameter $\beta$.

Consider the equation
\[
  \begin{split}
&\partial_t Z_\beta=\tfrac12\Delta Z_\beta+\beta Z_\beta V(t,x),\quad\mbox{
  where }\beta>0,\\
&
Z(0,x)\equiv1.
\end{split}
\]
Here $V(\cdot,\cdot)$ is a $1+1$-dimensional, spacetime white noise. By \cite[Corollary 2.5 ]{CD}, we know that 
\[
\EE[Z_\beta(t,0)^2]=2e^{\beta^4t/4}\int_{-\infty}^{\beta^2\sqrt{t/2}} \frac{1}{\sqrt{2\pi}}e^{-y^2/2}dy=:f(\beta,t).
\]
Let $\mathcal{N}_t$ be the Fourier spectrum   of $Z_\beta(t,0)$. It is a general fact that 
\begin{equation}\label{e.252}
\EE e^{-s \cN_t}=\frac{\EE Z_{\beta}(t,0) \cP_s Z_\beta(t,0) }{\EE Z_{\beta}(t,0)^2}, \quad\quad s\geq0,
\end{equation}
where $(\cP_s)_{s\geq0}$ is the Ornstein-Uhlenbeck semigroup
associated with the noise $V(\cdot,\cdot)$. By 
\cite[Page 55, Equation (1.67)]{nualart}, it is given by
\begin{equation}\label{050602-25}
\cP_s F(V)=\EE \left[F\Big(e^{-s}
V +\sqrt{1-e^{-2s}} \tilde V \Big)\big|V\right],\quad s\ge0,
\end{equation}
where $F(\cdot)$ is arbitrary square integrable function of $V$. 

To see why \eqref{e.252} holds, we note that, for random variable of the form $X_t=\sum_{n\geq0} I_n (f_{t,n})$ where $I_n(f_{t,n})$ is the $n-$th order chaos, given in \eqref{010402-25}, we have 
\[
\cP_s X_t=\sum_{n\geq0} e^{-ns} I_n (f_{t,n}),
\]
see \cite[Page 54, Definition 1.4.1]{nualart}. This directly implies \eqref{e.252}.

We claim that
\begin{equation}\label{e.claim}
\cP_s Z_\beta(t,0)=Z_{\beta e^{-s}}(t,0),
\end{equation}
which implies that 
\begin{equation}\label{e.253}
\EE Z_{\beta}(t,0) \cP_s Z_\beta(t,0)=\EE |\cP_{s/2} Z_\beta(t,0)|^2=\EE Z_{\beta e^{-s/2}}(t,0)^2.
\end{equation}
To prove the claim, we first assume $V$ is smooth in the spatial variable and by the Feynman-Kac representation the solution can be written as
\[
  Z_\beta(t,0)=\E\exp\left\{\beta\int_0^t V(r,x+B_{t-r})dr-\tfrac12\beta^2R(0)t\right\},
\]
where the expectation $\E$ is on the standard Brownian motion $B$
(independent of $V$) starting from the origin. By \eqref{050602-25}, we have 
\[
\mathcal{P}_s Z_\beta(t,0)=\EE\left[\E \exp\left\{\beta\int_0^t V_s(r,x+B_{t-r})dr-\tfrac12\beta^2R(0)t\right\}\bigg| V\right],
\]
where $V_s$ was defined in \eqref{e.defVs}. Interchange the order of expectation and because we have (for almost every realization of $B$) 
\[
  \begin{split}
&\EE\left[\exp\left\{\beta\int_0^t
    V_s(r,x+B_{t-r})dr-\tfrac12\beta^2R(0)t\right\}\bigg|V\right]\\
&
=\exp\left\{\beta
  e^{-s}\int_0^t V(r,x+B_{t-r})dr-\tfrac12\beta^2e^{-2s}R(0)t\right\},
\end{split}
\]
we conclude that $\mathcal{P}_s Z_\beta(t,0)=Z_{\beta e^{-s}}(t,0)$ 
when $V$ is smooth in the spatial variable. An approximation leads to the same conclusion when $V$ is a spacetime white noise, so \eqref{e.claim} is proved.

Combining \eqref{e.252} and \eqref{e.253} we have
\begin{equation}\label{e.laplaceSHE}
\begin{aligned}
\EE e^{-s \cN_t}=\frac{f(\beta e^{-s/2},t)}{f(\beta,t)}=\exp\left\{\frac{\beta^4t}{4}(e^{-2s}-1)\right\}\frac{\int_{-\infty}^{\beta^2 e^{-s}\sqrt{t/2}}\tfrac{1}{\sqrt{2\pi}}e^{-y^2/2}dy}{\int_{-\infty}^{\beta^2\sqrt{t/2}} \tfrac{1}{\sqrt{2\pi}}e^{-y^2/2}dy}.
\end{aligned}
\end{equation}

To prove the convergence in distribution, we use an analytic continuation argument. Define $g(z)=\EE e^{-z \mathcal{N}_t}$ for $z\in \mathbb{C}$ such that $\mathrm{Re}(z)\geq0$. Since $\mathcal{N}_t\geq0$, we know that $g$ is analytic in $\{z: \mathrm{Re}(z)>0\}$ and continuous up to $\mathrm{Re}(z)=0$. We also know that 
\[
g(z)=\exp\left\{\frac{\beta^4t}{4}(e^{-2z}-1)\right\}\frac{\int_{-\infty}^{\beta^2 e^{-z}\sqrt{t/2}}\tfrac{1}{\sqrt{2\pi}}e^{-y^2/2}dy}{\int_{-\infty}^{\beta^2\sqrt{t/2}} \tfrac{1}{\sqrt{2\pi}}e^{-y^2/2}dy}, \quad\quad \mbox{ for  } z=s\geq0.
\]
It is clear that the r.h.s. of the above display can be extended to an analytic function on $\mathbb{C}$: one write the numerator as 
\[
\tfrac12+\int_0^{\beta^2 e^{-z}\sqrt{t/2}}\tfrac{1}{\sqrt{2\pi}}e^{-y^2/2}dy
\] 
and the integration $dy$ is along any contour on $\mathbb{C}$
connecting $0$ and $\beta^2 e^{-z}\sqrt{t/2}$.
The above implies that 
\[
\EE e^{i\theta \mathcal{N}_t}= \exp\left\{\frac{\beta^4t}{4}(e^{-2i\theta
}-1)\right\}\frac{\tfrac12+\int_{0}^{\beta^2
    e^{-i\theta}\sqrt{t/2}}\tfrac{1}{\sqrt{2\pi}}e^{-y^2/2}dy}{\int_{-\infty}^{\beta^2\sqrt{t/2}}
  \tfrac{1}{\sqrt{2\pi}}e^{-y^2/2}dy},\quad \theta \in\R.
\]
After  centering and rescaling we conclude that for any $\theta \in\R $, 
\[
\EE \exp\left\{i\theta
  \frac{\cN_t-\tfrac{\beta^4t}{2}}{\sqrt{t}}\right\}\to
\exp\left\{-\frac{\beta^4\theta^2}{2}\right\},\quad \mbox{as }t\to+\infty.
\]
The proof is complete.

\begin{remark}
We expect similar results to hold for other initial data or covariance functions, in a certain strong disorder regime. The identity $\EE e^{-s\mathcal{N}_t}=\tfrac{\EE Z_{\beta e^{-s/2}}(t,0)^2}{\EE Z_{\beta}(t,0)^2}$ holds in all cases. By the second moment formula for the stochastic heat equation, we know that the growth of the second moment is related to the principle eigenvalue of the operator $\Delta+\beta^2 R(\cdot)$, denoted by $\Lambda(\beta)$. In other words, one may expect that $\EE e^{-s\mathcal{N}_t}\sim e^{\Lambda(\beta e^{-s/2})t-\Lambda(\beta)t}$, and, for small $s$,  a formal expansion leads to 
\[
\EE e^{-s\mathcal{N}_t}\sim e^{\Lambda(\beta e^{-s/2})t-\Lambda(\beta)t}\sim e^{\Lambda'(\beta)\beta(e^{-s/2}-1)t}.
\]
The last expression is precisely the Laplace transform of a Poisson distribution with parameter of order $t$.
\end{remark}



\begin{thebibliography}{99}


\bibitem{chatterjee}
Sourav Chatterjee. \emph{Superconcentration and related topics}, Vol. 15. Cham: Springer, 2014.

\bibitem{CD}
Le Chen and Robert C. Dalang. \emph{Moments and growth indices for the nonlinear stochastic heat equation with rough initial conditions}, Annals of Probability (2015): 3006-3051.



\bibitem{GH24}
Shirshendu Ganguly and Alan Hammond. \emph{Stability and chaos in dynamical last passage percolation}, Communications of the American Mathematical Society 4.09 (2024): 387-479.

\bibitem{GS14}
Christophe Garban and Jeffrey E. Steif. \emph{Noise sensitivity of Boolean functions and percolation}, Vol. 5. Cambridge University Press, 2014.

\bibitem{GK21}
Yu Gu  and Tomasz Komorowski. \emph{Gaussian fluctuations from random Schr\"odinger equation}, Communications in Partial Differential Equations 46.2 (2021): 201-232.

\bibitem{nualart}
David Nualart. \emph{Malliavin calculus and related topics}, (2006).

\end{thebibliography}
 \end{document}